\documentclass[12pt]{article}
\usepackage{amsmath,amssymb}
\usepackage{graphics}

\font\tenmsb=msbm10 \font\sevenmsb=msbm7 \font\fivemsb=msbm5

\catcode`\@=11 \ifx\amstexloaded@\relax\catcode`\@=\active
\endinput\else\let\amstexloaded@\relax\fi
\def\spaces@{\space\space\space\space\space}
\def\spaces@@{\spaces@\spaces@\spaces@\spaces@\spaces@}
\def\space@.{\futurelet\space@\relax}
\space@.
\def\Err@#1{\errhelp\defaulthelp@\errmessage{AmS-TeX error: #1}}
\def\relaxnext@{\let\next\relax}
\def\accentfam@{7}
\def\noaccents@{\def\accentfam@{0}}
\def\Cal{\relaxnext@\ifmmode\let\next\Cal@\else
\def\next{\Err@{Use \string\Cal\space only in math mode}}\fi\next}
\def\Cal@#1{{\Cal@@{#1}}}
\def\Cal@@#1{\noaccents@\fam\tw@#1}
\def\Bbb{\relaxnext@\ifmmode\let\next\Bbb@\else
\def\next{\Err@{Use \string\Bbb\space only in math mode}}\fi\next}
\def\Bbb@#1{{\Bbb@@{#1}}}
\def\Bbb@@#1{\noaccents@\fam\msbfam#1}
\newfam\msbfam
\textfont\msbfam=\tenmsb \scriptfont\msbfam=\sevenmsb
\scriptscriptfont\msbfam=\fivemsb

\def\beq{\begin{equation}}
\def\eeq{\end{equation}}

\def\qedbox{$\rlap{$\sqcap$}\sqcup$}

\newtheorem{thm}{Theorem}[section]
\newtheorem{lem}{Lemma}[section]

\newtheorem{prop}{Proposition}[section]

\newtheorem{rem}[thm]{Remark}

\@addtoreset{equation}{section}

\newcommand{\hess}{\mbox{Hess}}

\newcommand{\ric}{\mbox{Ric}}

\newcommand{\dv}{\mbox{div}}
\newcommand{\tr}{\mbox{trace}}
\newcommand{\vol}{\mbox{vol}}

\newcommand{\proof}{\par\medbreak\it Proof: \rm}

\def\proof{\noindent{\bf \it Proof.\ \ }}

\begin{document}

\bibliographystyle{plain}
\setlength{\baselineskip}{.51cm}
\title
{On Shrinking Gradient Ricci Solitons With Nonnegative Sectional Curvature}
\author{Mingliang Cai}
\date{}
\maketitle

\begin{abstract}
Perelman proved that an open $3$-dimensional shrinking gradient Ricci soliton with bounded nonnegative sectional curvature is a quotient of $S^2 \times \mathbb{R}$ or $\mathbb{R}^3$. We extend this result to higher dimensions with a decay condition on the Ricci tensor. 
\end{abstract}

\vskip 1.0cm

 \let\thefootnote\relax

 \section{Introduction} \[\] \vspace{-1.2cm}
 
A gradient Ricci soliton is a Riemannian manifold $(M, g)$ together with a smooth function $f$ such that
\[\ric+\hess f\,=\,\lambda \,g,\]
where $\lambda$ is a constant. It is called shrinking, steady and expanding when $\lambda >0$, $\lambda =0$ and $\lambda<0$ respectively.

Gradient Ricci solitons are self-similar solutions of Hamilton's Ricci flow and play a vital role in the analysis of singularities of the flow. In dimension 2, Hamilton \cite{H1} completely classified shrinking gradient Ricci solitons with bounded curvature and proved that they are the sphere, the projective space and the Euclidean space with constant curvature.
In dimension 3, Ivey \cite{I} proved that compact shrinking gradient Ricci solitons have positive sectional curvature and Perelman \cite{P} proved that shrinking gradient Ricci solitons with bounded nonnegative sectional curvature are quotients of $S^3$, $S^2 \times \mathbb{R}$ or $\mathbb{R}^3$.

In higher dimensions, there have been many results in the last several years. Chen \cite{Ch} showed that a complete shrinking gradient Ricci soliton has nonnegative scalar curvature. Ni and Wallace \cite{NW} gave the classification of shrinking gradient Ricci solitons with nonnegative Ricci curvature and zero Weyl tensor. Petersen and Wylie \cite{PW2} and independently,   Cao, Wang and Zhu \cite{CWZ},  classified the shrinking gradient Ricci solitons with zero Weyl tensor. Fern\'{a}ndez-L\'{o}pez and Garcia-Rio \cite{FG} considered solitons with harmonic Weyl tensor. In 
\cite{PW1}, several natural curvature conditions are given that characterize gradient Ricci solitons of the flat vector bundle $N \times_{\Gamma} \mathbb{R}^m$, where $N$ is an Einstein manifold, $\Gamma$ acts freely on $N$ and by orthogonal transformations on $\mathbb{R}^m$, and $f = \frac{1}{4} d^2$ with $d$ being the distance on the flat fiber to the base. In particular, it is shown in \cite{PW1} that a shrinking gradient Ricci soliton is rigid, i.e., of the form $N \times_{\Gamma} \mathbb{R}^m$, if the scalar curvature is constant and the sectional curvature of the plane containing $\nabla f$ is nonnegative.  As a consequence of a theorem of B\"{o}hm and Wilking (\cite{BW}), the gradient Ricci solitons with positive curvature operators are trivial. In view of this and the aforementioned result of Perelman, one naturally asks to what extend shrinking gradient Ricci solitons with nonnegative sectional curvature are rigid.  Our first result in this paper is the rigidity under a decay condition on $|D\ric|$, extending Perelman's result to higher dimensions.  In all theorems we scale the metric so that  $\lambda = \frac{1}{2}$.

\vspace{.2cm}

\begin{thm}\label{main0.5}\,\, Let $(M, g, f)$ be a complete non-compact shrinking gradient Ricci soliton with bounded nonnegative sectional curvature.  Assume that there exists $\delta > 0$ such that 
 \[ \int_{M}e^{\delta f} |D \ric| \,d\vol _g\,< \infty.\] Then $(M^n, g)$ is isometric to $N \times_{\Gamma} \mathbb{R}^m$, where $N$ is a compact Einstein manifold. 
\end{thm}

This is, to our knowledge, the first rigidity result in high dimensions without assumptions on the Weyl tensor. The potential function $f$ is known to grow quadratically with respect the distance from a fixed point, so our condition on $ D\ric$ says that it decays exponentially. Our proof can be seen to work  under the assumption that $D\ric$ decays polynomially with a degree depending on other geometric quantities. 

The Cheeger-Gromoll Soul Theorem states that an open manifold with nonnegative sectional curvature is diffeomorphic to a vector bundle over a compact submanifold called a soul. The pull-back metric on the bundle can be highly twisted. However, if there exists a gradient soliton structure on such a bundle, then,  by Theorem \ref{main0.5},  the metric has to be locally trivial, provided that the decay condition is satisfied. The decay condition on $D\ric$ in Theorem \ref{main0.5} is imposed in the region where $f$ is large. Our next result deals with the rigidity under a condition on $D\ric$ imposed in the region where $f$ is small. 

\begin{thm}\label{main1}  \,\, Let $(M^n, g, f)$ be a complete shrinking gradient Ricci soliton with bounded nonnegative sectional curvature. Assume that the minima of $f$ is a smooth compact non-degenerate critical submanifold, $D\ric$ and $D^2\ric$ vanish on the minima, then $(M^n, g)$ is non-compact and isometric to $N \times_{\Gamma} \mathbb{R}^m$, where $N$ is a compact Einstein manifold.
\end{thm}

\vspace{.2cm}

We derive some basic formulas in section 2, and prove theorems \ref{main0.5} and \ref{main1} in sections 2 and 3 respectively. \vspace{.2cm}

{\bf Acknowledgement.} I thank Professors Peter Petersen and DaGang Yang for their interests in this work and for helpful discussions. I thank Professor Ovidiu Munteanu for pointing out an error in an earlier version of the paper.  I also thank the referee for the thorough review and helpful suggestions.

\section{Basic Formulas}\[\] \vspace{-1.2cm}

There are different conventions for the curvature tensor in the literature, to avoid the confusion, we state ours as follows.
The $(3,1)$ tensor $Rm(X,Y,Z)=Rm(X,Y)Z$ is defined as
\[Rm(X,Y)Z=D_X\,D_Y\,Z-D_Y\,D_X\,X-D_{[X,Y]}Z\]
and the $(4,0)$ tensor as
\[Rm(X,Y,Z,W)=<Rm(X,Y)Z,W>.\]
We use $\ric$ to denote the Ricci tensor and $R$ the scalar curvature. For a tangent vector $X$ at $p$, we use $\ric(X)$ to denote the vector  such that 
 \[<\ric(X),Y>=\ric(X,Y)\]
 for any vector  $Y$ at $p$.  For any smooth vector field $V$ and any smooth function $\phi$ on manifold $M$,  by $V(\phi)$, we mean $V(\phi) = d \phi (V) = < V, \nabla \phi>$.  
 In the remaining of the paper, we will rescale the metric and assume that our gradient Ricci soliton satisfies
    \[\ric+\hess f\,=\,\frac{1}{2}\,g.\]  
Since the curvature of $(M, g)$ is assumed to be bounded, there exists a flow $\Phi_t: M \rightarrow M$ defined for all time with $\Phi_0=\mbox{Id}$ and $\frac{\partial \Phi}{\partial t}\,=\nabla f$ ( p. 207, ~\cite{MT}). For $ t \in (\infty, 0)$, define $G(t)=|t|\,\Phi^{*}_{-\ln |t|}\,g$. Then $G(-1)=g$ and $G(t)$ satisfies 
 \[\ric(G(t)) + \hess f = \frac{1}{2\tau}\,G(t),\]
 where $\hess$ is taken with respect to the metric $G(t)$ and $\tau=|t|=-t$.

In the next lemma, we collect some well-known formulae.
  
 \begin{lem} \label{basic}\,\, On (M, G(t)), we have
 \begin{eqnarray*}
\mbox{(1)} \hspace{1cm} && dR=2\,\ric(\nabla f, \cdot)\\
\mbox{(2)}\hspace{1cm}&& |\nabla f|^2=\frac{f}{\tau} -R + \mbox{constant}\\ 
\mbox{(3)}\hspace{1cm}&&  \frac{R}{\tau} + <\nabla f, \nabla R>\,=\,\Delta R\,+\, 2\,|\ric|^2\\
\mbox{(4)}\hspace{1cm}&& \dv Rm (X,Y,Z)=Rm(\nabla f, X, Y, Z)\\
\mbox{(5)}\hspace{1cm}&& D_{X}\ric(Y,Z)=D_{Y}\ric(X,Z)-Rm(X,Y,\nabla f, Z),
\end{eqnarray*}
where $\dv Rm (X,Y,Z)= \tr_{1,2}DRm(\cdot,\cdot, X,Y,Z)$.
\end{lem}
  \vspace{.2in}

 \proof  The derivation of (1)-(3) can be found in \cite{H2} and (4)-(5) in \cite{PW2}.
 
  \hfill{} \qedbox

 \vspace{.2in}

 \begin{lem} \label{delta-ricci-square} \,\, On $(M, g)$, the following holds. 
  \begin{eqnarray*}
  \Delta |\ric|^2 = 2 |D\ric|^2 + 2 |\ric|^2 + \nabla f (|\ric|^2) - 4 K_{ij}\lambda_i\lambda_j,
  \end{eqnarray*}
  where $\lambda_i$ are the eigenvalues of the Ricci tensor and $K_{ij}$ is the sectional curvature of the plane spanned by the eigenvectors belonging to $\lambda_i$ and $\lambda_j$ respectively.
  \end{lem}
   \vspace{.2in}
   
   \proof  This follows from the following formula derived in Lemma 2.1 in \cite{PW2}.
   \[ \Delta \ric = D_{\nabla f} \ric + \ric -2 \sum_{k=1}^n\, Rm(\cdot, e_k, \ric(e_k), \cdot).\]
      
       \hfill{} \qedbox

 
  \vspace{.2in}

 Throughout the computations in the paper, we assume  $\{e_1, ..., e_{n}\}$ is an orthonormal basis 
  in a neighborhood of a fixed point $x$ with $D_{e_i}e_{j}\,(x)=0$ and further assume that each $e_i$ is an eigenvector 
  of $\ric$ at $x$ corresponding to the eigenvalue $\lambda_i$. Such a basis always exists. We also use the Einstein summation convention (unless otherwise specified).

 \begin{lem} \label{div-ric-nabla-r-1} \,\, On $(M, g)$, we have
  \begin{eqnarray*}
  \dv(\ric(\nabla R))=\nabla f (|\ric|^2)+\frac{1}{2}|\nabla R|^2 -2<Z, \nabla f> + |\ric|^2 - 2\sum_{i}\,\lambda_i ^3 ,
 \end{eqnarray*}
 where $Z=\ric(e_i,e_j)Rm (\nabla f, e_i, e_j)$.
 \end{lem}

 \proof  The following computations are done at $x$. We have from Lemma \ref{basic} 
 \begin{eqnarray*}
 D_{e_i}\ric (\nabla R, e_i) &=& D_{\nabla R }\ric (e_i, e_i) - Rm (e_i, \nabla R, \nabla f, e_i)\\
 &=& |\nabla R|^2 - \ric(\nabla R, \nabla f) =  \frac{1}{2}|\nabla R|^2.
 \end{eqnarray*} 
 We then obtain  \vspace{-0.1cm} 
 \begin{eqnarray*}
 \dv(\ric(\nabla R))&=&<D_{e_i}\ric(\nabla R), e_{i}>=e_{i}\ric(\nabla R, e_i)\\ 
  &=&D_{e_i}\ric (\nabla R, e_i)+\ric(D_{e_i}\nabla R, e_i)\\
 &=& \frac{1}{2}|\nabla R|^2 + \ric(e_i, e_j)<D_{e_i}\nabla R, e_j>\\
 &=& \frac{1}{2}|\nabla R|^2 + 2\ric(e_i, e_j)<D_{e_i}\ric(\nabla f), e_j>\\
 &=& \frac{1}{2}|\nabla R|^2 + 2\ric(e_i, e_j)e_i \ric(\nabla f, e_j)\\
 &=& \frac{1}{2}|\nabla R|^2 + 2\ric(e_i, e_j)[D_{e_i}\ric(\nabla f, e_j)+\ric(D_{e_i}\nabla f, e_j)].
 \end{eqnarray*}
That is,
 \begin{eqnarray}
 \dv(\ric(\nabla R))= \frac{1}{2}|\nabla R|^2 + 2\ric(e_i, e_j)[D_{e_i}\ric(\nabla f, e_j)+\ric(D_{e_i}\nabla f, e_j)].  \label{eq:bu4}
 \end{eqnarray} 
From the soliton equation \[ \ric+\hess f\,=\,\frac{1}{2}\,g\]  
it follows that 
\[ D_{e_i}\nabla f = \frac{1}{2} e_i - \ric(e_i) = \frac{1}{2} e_i - \lambda_i e_i, \]
where we have used the assumption that $e_i$ is an eigenvector of $\ric$ at $x$ belonging to the eigenvalue $\lambda_i$. 
Hence, 
\begin{eqnarray}
2\ric(e_i, e_j)\ric(D_{e_i}\nabla f, e_j) = 2 (\frac{1}{2}-\lambda_i) [\ric(e_i,e_j)]^2 = 2\lambda_i^2 (\frac{1}{2} - \lambda_i).  \label{eq:bu5}
\end{eqnarray}
Lemma 2.1 (5) implies that 
\[D_{e_i}\ric(\nabla f, e_j) = D_{\nabla f}\ric(e_i, e_j) - Rm(e_i, \nabla f, \nabla f, e_j).\]
It follows that
\begin{eqnarray}
2\ric(e_i, e_j) D_{e_i}\ric(\nabla f, e_j)  &=& 2\ric(e_i, e_j) [D_{\nabla f}\ric(e_i, e_j) - Rm(e_i, \nabla f, \nabla f, e_j)] \nonumber \\
&=& 2\ric(e_i, e_j) D_{\nabla f}\ric(e_i, e_j) - 2< Z, \nabla f > \nonumber\\
&=& \nabla f (|\ric|^2) - 2 <Z, \nabla f>.     \label{eq:bu6}
\end{eqnarray}
Combining (\ref{eq:bu5}) and (\ref{eq:bu6}), we obtain that 
\begin{eqnarray*}
&&2\ric(e_i, e_j)[D_{e_i}\ric(\nabla f, e_j)+\ric(D_{e_i}\nabla f, e_j)] \\
 &=& \nabla f (|\ric|^2)-2<Z, \nabla f> +2 \lambda_i^2(\frac{1}{2}-\lambda_i).
 \end{eqnarray*}
 Substituting the above into (\ref{eq:bu4}) gives
 \begin{eqnarray*}
 \dv(\ric(\nabla R))
 &=&\frac{1}{2}|\nabla R|^2 + \nabla f (|\ric|^2)-2<Z, \nabla f> +2 \lambda_i^2(\frac{1}{2}-\lambda_i)\\
 &=& \frac{1}{2}|\nabla R|^2 + \nabla f (|\ric|^2)-2<Z, \nabla f> + |\ric|^2 - 2\sum_{i}\,\lambda_i ^3.  
 \end{eqnarray*}
 Lemma \ref{div-ric-nabla-r-1} is thus proved. 
 \hfill{} \qedbox

\begin{rem} \label{z-f} 
 $< Z, \nabla f> \geq 0$ , when the sectional curvature of $(M,g)$ is nonnegative. In fact, at $x$, $ <Z, \nabla f> = \lambda_i Rm (\nabla f, e_i, e_i, \nabla f)$.
 \end{rem}

The next lemma is a slight variation of Lemma \ref{div-ric-nabla-r-1}.

 \begin{lem} \label{div-ric-nabla-r-2} \,\, On $(M, g)$, we have
  \begin{eqnarray*}
\nabla f (|\ric|^2) = 2[< Z, \nabla f > + \sum_{i=1}^{n}\,\lambda_i (\lambda_i - \frac{1}{2})^2] + \frac{1}{2} <\nabla f, \nabla R > - \frac{1}{2} |\nabla R|^2 - \dv(D_{\nabla R} \nabla f).
 \end{eqnarray*}
  \end{lem}
  \vspace{.2in}

\proof  It follows from Lemma  \ref{div-ric-nabla-r-1} that 
\[\dv(\ric(\nabla R)) = \frac{1}{2}|\nabla R|^2 + \nabla f (|\ric|^2)-2<Z, \nabla f> + |\ric|^2 - 2\sum_{i}\,\lambda_i ^3.  \]
Using $\ric(\nabla R) = \frac{1}{2} \nabla R - D_{\nabla R}\nabla f$ and Lemma \ref{basic} (3), we have
\begin{eqnarray*}
&&\nabla f (|\ric|^2) = \frac{R}{2}- 2|\ric|^2+ 2\sum_{i}\,\lambda_i ^3\\
&& \hspace{2.5cm} + 2< Z, \nabla f >  + \frac{1}{2} <\nabla f, \nabla R > - \frac{1}{2} |\nabla R|^2 - \dv(D_{\nabla R} \nabla f).
\end{eqnarray*}
The lemma now follows as $\frac{R}{2}- 2|\ric|^2+ 2\sum_{i}\,\lambda_i ^3 = 2 \sum_{i=1}^{n}\,\lambda_i (\lambda_i - \frac{1}{2})^2.$

 \hfill{} \qedbox

 Combining Lemma \ref{div-ric-nabla-r-1}  with \ref{delta-ricci-square}  gives the following proposition.  \vspace{.2in}

  \begin{prop} \label{P} \,\, On $(M, g)$, 
   \begin{eqnarray*}
P=\frac{1}{2}\nabla f (|\ric|^2)+\frac{1}{2}|\nabla R|^2+\dv [\frac{1}{2}\nabla |\ric|^2-\ric(\nabla R)],
\end{eqnarray*}
where $P=K_{ij}(\lambda_i-\lambda_j)^2+|D\ric|^2+2<Z,\nabla f>$.
\end{prop}

  \vspace{.2in}

 \proof  Lemma  \ref{delta-ricci-square} implies that 
 \begin{eqnarray*}
 - 2 K_{ij}\lambda_i\lambda_j + |D\ric|^2  =   - \frac{1}{2}  \nabla f (|\ric|^2) -  |\ric|^2+ \dv(\frac{1}{2}\,\nabla |\ric|^2),
  \end{eqnarray*}
 while Lemma \ref{div-ric-nabla-r-1} implies that 
 \begin{eqnarray*}
  2\sum_{i}\,\lambda_i ^3 + 2<Z, \nabla f>  = \nabla f (|\ric|^2) +  |\ric|^2+ \frac{1}{2}|\nabla R|^2 -  \dv(\ric(\nabla R)).
 \end{eqnarray*}
 
 Adding the corresponding sides of the last two equations and noting that $2 \sum_{i} \, \lambda _i^3 - 2 \sum_{i,j}\, K_{ij}\lambda_i\lambda_j =  \sum_{i,j}\, K_{ij}(\lambda_i-\lambda_j)^2$, we obtain Proposition \ref{P}.
 
  \hfill{} \qedbox
 
   \vspace{.2in}
  
  \begin{rem} \label{P+} 
 Clearly, $P \geq 0$ , when the sectional curvature of $(M,g)$ is nonnegative.
 \end{rem}

The proof of Theorems \ref{main0.5} will use an alternative form of Proposition \ref{P} in which the term $|D\ric|^2$ is replaced by $|\dv Rm|^2$.  An integral from of next lemma is proved in 
\cite{Ca}.
  
 \begin{lem}\label{d-ric-dv-rm} On $(M, g)$, 
   \begin{eqnarray*}
|D\ric|^2 = |\dv Rm|^2 + 2 <Z, \nabla f> - \frac{1}{2} \nabla f (\ric|^2) + \dv( \frac{1}{2} \nabla |\ric|^2 - 2 Z).
\end{eqnarray*}
\end{lem}

\vspace{.2in}

\proof   \,\, As before, we fix 
 an orthonormal 
 basis, $\{e_1, ..., e_{n}\}$, in an neighborhood of a fixed point $x$ and assume that $D_{e_i}e_{j}\,(x)=0$ and that each $e_i$ is an eigenvector 
  of $\ric$ at $x$ corresponding to the eigenvalue $\lambda_i$. Recall that $Z = \ric(e_i,e_j)Rm (\nabla f, e_i, e_j)$, so at $x$,
 \begin{eqnarray*}
\dv (Z) &= &<D_{e_k} Z, e_k> = <D_{e_k} [\ric(e_i, e_j) Rm(\nabla f, e_i, e_j)] , e_k>\\
&=& e_k [\ric(e_i,e_j)] Rm(\nabla f, e_i, e_j, e_k) + \ric(e_i,e_j) <D_{e_k} [Rm(\nabla f, e_i, e_j)] , e_k>\\
&=& D_{e_k}\ric (e_i, e_j) \, Rm (\nabla f, e_i, e_j, e_k) + \ric(e_i, e_j) e_k [Rm(\nabla f, e_i, e_j , e_k)]\\
&=& D_{e_k}\ric(e_i, e_j) \dv Rm (e_i, e_j, e_k) + \ric(e_i, e_j) [D_{e_k} Rm (\nabla f, e_i, e_j , e_k) \\
&&\hspace{7.5cm}+ Rm (D_{e_k}\nabla f, e_i, e_j, e_k)]\\
&=& [D_{e_i}\ric(e_j, e_k) - Rm(e_k, e_i, \nabla f, e_j)] \dv Rm (e_i, e_j, e_k)  \\
&&\hspace{0.5cm}+ \ric(e_i, e_j) \dv Rm(e_j, e_i, \nabla f ) + \lambda_i\, Rm ((\frac{1}{2}-\lambda_k)e_k, e_i, e_i, e_k)\\
&=& [D_{e_i}\ric(e_j, e_k) \dv Rm (e_i, e_j, e_k)  + \dv Rm (e_j, e_i, e_k)\,\dv Rm (e_i, e_j, e_k)\\
&&\hspace{0.5cm}+ \ric(e_i, e_j) Rm(\nabla f, e_j, e_i, \nabla f ) + K_{ij} \lambda_i\, (\frac{1}{2}-\lambda_j).
  \end{eqnarray*}
In the above calculation, we have repeatedly used Lemma \ref{basic}. The lemma now follows from Lemma \ref{delta-ricci-square} and the following two identities whose 
proofs are easy.
\[D_{e_i}\ric(e_j, e_k) \dv Rm (e_i, e_j, e_k)=0\] and
\[ \dv Rm (e_j, e_i, e_k)\,\dv Rm (e_i, e_j, e_k)=\frac{1}{2} |\dv Rm|^2. \]

 \hfill{} \qedbox

  \vspace{.2in}

 Lemma \ref{d-ric-dv-rm}, together with Proposition \ref{P}, implies the following
 
  \begin{lem} \label{pre-Q} \,\, On $(M, g)$, 
   \begin{eqnarray*}
Q=\nabla f (|\ric|^2)+\frac{1}{2}|\nabla R|^2+\dv [2 Z -\ric(\nabla R)],
\end{eqnarray*}
where $Q=K_{ij}(\lambda_i-\lambda_j)^2+|\dv Rm|^2+4<Z,\nabla f>$.
\end{lem}

 \hfill{} \qedbox

  \vspace{.2in}
  
   \begin{rem} \label{Q+} 
 We note that $Q \geq 0$ , when the sectional curvature of $(M,g)$ is nonnegative.
 \end{rem}

  The next lemma deals with the term $\nabla f (|\ric|^2)$ in Lemma \ref{pre-Q}.

 \begin{lem} \label{nabla-f-ric}\,\,  On $(M, g)$, 
   \begin{eqnarray}
&& \nabla f (|\ric|^2)=\frac{1}{2} |\nabla R |^2+\frac{1}{2}<\nabla f, \nabla R>+\frac{1}{2} \nabla f (<\nabla f, \nabla R>)  \nonumber\\
&&\hspace{5cm}+\dv[D_{\nabla R}\,\nabla f\,-\, \frac{1}{2}  \nabla <\nabla f, \nabla R>)].  \label{eq:lem-nabla}
\end{eqnarray}
\end{lem}
\vspace{.2in}

\proof  It follows from Lemma \ref{basic}  (3) and (1) that
   \begin{eqnarray*}
\frac{1}{2}\nabla f(\Delta R) &=& - \nabla f (|\ric|^2) + \frac{1}{2} <\nabla f, \nabla R> + \frac{1}{2} \nabla f (<\nabla f, \nabla R>).\\
\end{eqnarray*}
The Bochner-Weitzenb\"{o}ck formula implies that
   \begin{eqnarray*}
\dv[  \frac{1}{2}  \nabla <\nabla f, \nabla R> ]&=&\frac{1}{2}\Delta <\nabla f, \nabla R>\\
&=& <\hess f, \hess R>+\frac{1}{2}\nabla f (\Delta R)+\frac{1}{2}\nabla R (\Delta f) + \ric(\nabla f, \nabla R)\\
&=&<\hess f, \hess R>+\frac{1}{2}\nabla f (\Delta R)+\frac{1}{2}\nabla R (\frac{n}{2}-R) +\frac{1}{2}|\nabla R|^2 \\
&=&<\hess f, \hess R>+\frac{1}{2}\nabla f (\Delta R).
\end{eqnarray*}
But,
 \begin{eqnarray*}
\dv(D_{\nabla R}\nabla f)&=&<D_{e_i}D_{\nabla R}\nabla f, e_i>=e_i<D_{\nabla R}\nabla f, e_i>=e_i<D_{e_i}\nabla f, \nabla R>\\
&=& <D_{e_i}(\frac{1}{2}e_i-\ric(e_i)), \nabla R>+<\hess f, \hess R>\\
&=& - D_{e_i}\ric(e_i, \nabla R) + <\hess f, \hess R>\\
&=&-\frac{1}{2}|\nabla R|^2 + <\hess f, \hess R>.
\end{eqnarray*}
The lemma follows. 

 \hfill{} \qedbox
\vspace{.2in}

 We now have the following proposition which will be used in the proof of Theorems \ref{main0.5}.
 
     \begin{prop} \label{Q}\,\,  On $(M, g)$, 
   \begin{eqnarray*}
   Q &= & |\nabla R|^2 +\frac{1}{2} <\nabla f, \nabla R> + \frac{1}{2} \nabla f[ <\nabla f, \nabla R>]\\
&&\hspace{3cm}+\dv[2 Z -\ric(\nabla R)+ D_{\nabla R}\,\nabla f\,-\,\frac{1}{2} \nabla  <\nabla f, \nabla R> ].
\end{eqnarray*}
\end{prop}

\proof  This is merely a consequence of Lemmas \ref{pre-Q} and \ref{nabla-f-ric}.  

 \hfill{} \qedbox

\vspace{.2in}

\section{Proof of Theorem \ref{main0.5}} \[\] \vspace{-1.2cm}

We will use $\phi$ to denote a real-valued nonnegative $C^4$ function on $\mathbb{R}$ and write  $\phi\circ f$ as $\phi(f)$. We will show that $R$ is a constant function and then appeal to \cite{PW1} to complete the proof.  We begin with the following proposition.

\begin{prop} \label{for-thm1.1} On $(M, g)$, 
  \begin{eqnarray}
&& \phi(f)\, Q = \frac{1}{2} <\nabla f, \nabla R> [(\phi - \phi^{\prime}) (f) - ( \phi + \phi^{\prime})(f) \Delta f  \nonumber \\
&&\hspace{7.8cm} -  (\phi^{\prime\prime} + \phi^{\prime})(f) |\nabla f |^2  ] \nonumber \\
 &&\hspace{2.3cm}  +  (\phi + \phi^{\prime}) (f) |\nabla R|^2 - 2 \phi^{\prime} < Z, \nabla f >  +  \dv (X),  \label{eq:for-thm1.1.1}
  \end{eqnarray} 
where 
\begin{eqnarray*}
X &=&  \frac{1}{2} <\nabla f, \nabla R> (\phi^{\prime} + \phi)(f) \nabla f  \\
&& \hspace{.2in} +  \phi(f) \,[2 Z -\ric(\nabla R)+ D_{\nabla R}\,\nabla f\,-\,\frac{1}{2} \nabla  <\nabla f, \nabla R>].
\end{eqnarray*}
\end{prop}

\proof  We  multiply each side of the equation in Proposition \ref{Q} by  $\phi(f)$ to get 
\begin{eqnarray*}
 && \phi(f) Q = \phi(f)  |\nabla R|^2 + \frac{\phi(f)}{2} <\nabla f, \nabla R> + \frac{\phi(f)}{2} \nabla f [ <\nabla f, \nabla R>]   \\
 &&\hspace{1.6cm} -\phi^{\prime}(f)\, < 2 Z -\ric(\nabla R)+ D_{\nabla R}\,\nabla f\,-\,\frac{1}{2} \nabla  <\nabla f, \nabla R>\,,\, \nabla f>   \\
  &&\hspace{2.4cm} + \dv \{\phi(f) \,[2 Z -\ric(\nabla R)+ D_{\nabla R}\,\nabla f\,-\,\frac{1}{2} \nabla  <\nabla f, \nabla R>]\} . 
  \end{eqnarray*}
It follows from the soliton equation and Lemma \ref{basic} (1) that
\begin{eqnarray*}
< -\ric(\nabla R)+ D_{\nabla R}\,\nabla f\, , \, \nabla f >  &=& < \frac{1}{2} \nabla R - 2 \ric(\nabla R) , \nabla f>\\ 
&=& \frac{1}{2} < \nabla f, \nabla R> - |\nabla R|^2.
\end{eqnarray*}
We thus obtain 
 \begin{eqnarray}
&&  \phi(f) Q  = (\phi+\phi^{\prime}) (f) | \nabla R|^2+\frac{\phi-\phi^{\prime}}{2} (f) <\nabla f, \nabla R> - 2 \phi^{\prime} < Z, \nabla f > \nonumber \\
&& \hspace{6.5cm} +\frac{ \phi + \phi^{\prime}}{2} (f) \nabla f (<\nabla f, \nabla R>)  \nonumber\\
&& \hspace{1cm} +\dv \{ \phi(f) \,[2 Z -\ric(\nabla R)+ D_{\nabla R}\,\nabla f\,-\,\frac{1}{2} \nabla  <\nabla f, \nabla R>] \}. \label{eq:for-thm1.1.2}
\end{eqnarray}
Now, we observe that 
\begin{eqnarray*}
&& (\phi+\phi^{\prime})(f) \nabla f (<\nabla f, \nabla R>) \\
&& \hspace{1cm} = <\nabla <\nabla f, \nabla R> , (\phi^{\prime} + \phi)(f)  \nabla f >\\
&& \hspace{1cm} = - <\nabla f, \nabla R> [ (\phi^{\prime} + \phi)(f) \Delta f +  (\phi^{\prime\prime} + \phi^{\prime})(f) |\nabla f |^2 ] \\
&& \hspace{5cm} + \dv [ <\nabla f, \nabla R> (\phi^{\prime} + \phi)(f)\nabla f ].
 \end{eqnarray*}
Substituting the above into (\ref{eq:for-thm1.1.2}), we obtain (\ref{eq:for-thm1.1.1}). Proposition \ref{for-thm1.1} is thus proved.

 \hfill{} \qedbox
\vspace{.2in}

The idea now is to choose an appropriate function $\phi$ and integrate (\ref{eq:for-thm1.1.1}) over $M$. The divergence term, after integration, vanishes because of the fall-off condition we impose. The right hand side will then be nonpositive while the left is always nonnegative, and consequently, $R$ is a constant. Theorem \ref{main0.5} follows from \cite{PW1}. 
\vspace{.2in}

{\em Proof of Theorem \ref{main0.5}} \,\,\,\,  We normalize $f$ by adding a constant so that Lemma \ref{basic} (2) takes the form $| \nabla f |^2 = f -R$. Since $R \geq 0$, we always have $ |\nabla f |^2 \leq f$. On the other hand, since $R$ is assumed to be bounded and  $f$ grows quadratically with respect to the distance from a fixed point (\cite{CZ}, \cite{Nb}), we have $ |\nabla f|^2 \geq \frac{1}{2}f$, when $f$ is sufficiently large. Thus, there exists $T > 2 $ so that when $f \geq T$, 
\begin{eqnarray}
 \frac{1}{2} f  \leq |\nabla f|^2 \leq f.    \label{eq:f-nabla-f}
 \end{eqnarray}
Fix $ 0 < \eta < \delta$  and define $\phi: \mathbb{R} \rightarrow \mathbb{R}$ by $\phi(t) = 0 $ for $ t \leq T$, and $\phi(t) = ( t -T)^{k} e^{\eta t}$ for $ t \geq T$, where $k$ is a sufficiently large number to be determined. Throughout this section, we will use this $\phi$ in (\ref{eq:for-thm1.1.1}). By our fall-off assumption, there exists a sequence $t_i \rightarrow \infty$ such that
\[ \int_{f=t_i} \, e^{\delta f}\, \frac{1}{ |\nabla f|}\, |D\ric|\, \rightarrow 0, \,\,\,\, \mbox{as}\,\, i \rightarrow \infty.\]  From this, we now deduce that 
\begin{eqnarray}
 \int_{ f \leq t_i}\, \dv (X)  \, =\, \int_{f=t_i}\, \frac{ <X, \nabla f >}{ |\nabla f|} \, \rightarrow  0, \,\,\,\, \mbox{as}\,\,  i \rightarrow \infty.   \label{eq:div-0}
 \end{eqnarray}
To this end, we look at each of the five terms in $X$ and denote by $X_i$ the $i^{\mbox{th}}$ term. Then, when $f > T$, 
\begin{eqnarray*}
\frac{ | <X_1, \nabla f>| }{|\nabla f|} = \frac{1}{2} |< \nabla f, \nabla R> | (\phi^{\prime}+\phi)(f) |\nabla f|  \leq  C_1\, f ^{k+1} e^{\eta f} |\nabla R|,
\end{eqnarray*}
where $C_1$ is a constant depending only on $k$ and $\eta$. Now by the Cauchy-Schwarz inequality, 
\begin{eqnarray*}
| D\ric |^2 = \sum_{i,j,k}\, [ D_{e_i} \ric (e_j, e_k) ]^2 \geq \frac{1}{n} \sum_i [ \sum_{j} D_{e_i} \ric ( e_j, e_j) ]^2 = \frac{1}{n} |\nabla R|^2.
\end{eqnarray*} Thus, \[ |\nabla R| \leq \sqrt{n} | D\ric|. \]
Hence, 
\begin{eqnarray*}
\frac{ |<X_1, \nabla f>| }{|\nabla f|}   \leq  C_1\, \sqrt{n}\, f ^{k+1} e^{\eta f} |D\ric|. 
\end{eqnarray*}
Integrating the above over $\{f = t_i\}$ and noting that \[C_1\, \sqrt{n}\, f^{k+1} e^{\eta f} |D\ric| \, \leq e^{\delta f}\frac{ |D \ric | }{ |\nabla f|},\] when $f$ is sufficiently large, we conclude
\begin{eqnarray*}
\int_{f=t_i}\, \frac{ | <X_1, \nabla f>| }{|\nabla f|}  \, \rightarrow  0, \,\,\,\, \mbox{as}\,\,  i \rightarrow \infty.
\end{eqnarray*}
Now note that  $ <X_2, \nabla f > = 2 \phi < Z, \nabla f > = 2 \phi \sum_{i}\, \lambda_i\, Rm(\nabla f, e_i, e_i, \nabla f)$. Since $\ric$ is assumed to be bounded and since the sectional curvature is nonnegative, 
\begin{eqnarray*}
\frac{ |<X_2, \nabla f>| }{|\nabla f|}   \leq  C_2\, \, f^{k-\frac{1}{2}} e^{\eta f} \ric(\nabla f, \nabla f) = C_2\, \, f^{k-\frac{1}{2}} e^{\eta f} \frac{1}{2} < \nabla f, \nabla R>, 
\end{eqnarray*}
where $C_2$ is a constant dependent only on the bound of $\ric$ and the last equality follows from Lemma \ref{basic}. Hence, when $f$ is sufficiently large, 
\begin{eqnarray*}
\frac{ |<X_2, \nabla f>| }{|\nabla f|}   \leq  \frac{1}{2}\,  C_2\, \, f^{k} e^{\eta f} | \nabla R |  \, \leq e^{\delta f}\frac{ |D \ric | }{ |\nabla f|}.
\end{eqnarray*}
It then follows that
\begin{eqnarray*}
\int_{f=t_i}\, \frac{ | <X_2, \nabla f>| }{|\nabla f|}  \, \rightarrow  0, \,\,\,\, \mbox{as}\,\,  i \rightarrow \infty.
\end{eqnarray*}
The arguments for other $X_i$'s are similar, we will skip $X_3$ and $X_4$.  Now look at $X_5$. Repeatedly using Lemma \ref{basic}(2), we see that 
\begin{eqnarray*}
< X_5, \nabla f > &=& -\frac{1}{2} \phi \nabla f (< \nabla f, \nabla R>) = - \phi \nabla f [\ric (\nabla f, \nabla f)] \\
&=& - \phi [\, D_{\nabla f} \ric (\nabla f, \nabla f) + 2 \ric (D_{\nabla f} \,\nabla f , \nabla f)]\\
&-& - \phi [\, D_{\nabla f} \ric (\nabla f, \nabla f) +  \ric ( \nabla f - \nabla R, \nabla f)]\\
&=& -\phi [ \, D_{\nabla f} \ric (\nabla f, \nabla f) +  \frac{1}{2} <\nabla f, \nabla R> - \ric ( \nabla R, \nabla f)].
\end{eqnarray*}
Since $|\nabla R|$ can be bounded by $|D\ric|$, we have $  | <X_5, \nabla f > | \leq C_5\, e^{\eta f} f^{k+3} | D \ric| $.  (\ref{eq:div-0}) then follows. 

\noindent To simplify notations, we put 
\begin{eqnarray*}
&& F  = \frac{1}{2} <\nabla f, \nabla R> [(\phi - \phi^{\prime}) (f) - ( \phi + \phi^{\prime})(f) \Delta f  \nonumber \\
&&\hspace{5.3cm} -  (\phi^{\prime\prime} + \phi^{\prime})(f) |\nabla f |^2  ] \nonumber \\
 &&\hspace{2.5cm}  +  (\phi + \phi^{\prime}) (f) |\nabla R|^2 - 2 \phi^{\prime} < Z, \nabla f > . 
 \end{eqnarray*}
Then,
\begin{eqnarray*}
\phi (f) Q  = F + \dv (X).
 \end{eqnarray*}
It follows easily from the arguments in the proof of (\ref{eq:div-0}) that $\int_{M}\, F \, d\vol_g\, < \, \infty$. We thus have
\begin{eqnarray}
\int_{M}\, \phi (f) Q  =  \int_{M} \, F.   \label{eq:bu1}
\end{eqnarray}
We now show that $\int_{M}\, F \, d\vol_g\, \leq 0$. First, we note that $- \Delta f  = R - \frac{n}{2}  \leq \Lambda$, where $\Lambda$ is an upper bound of $R$, hence 
$ - ( \phi + \phi^{\prime})(f) \Delta f  \leq \Lambda (\phi + \phi^{\prime})$, as $\phi$ and $\phi^{\prime}$ are both nonnegative. Next, we observe that, by Lemma \ref{basic}, \[ |\nabla R|^2 = 2 \ric ( \nabla f, \nabla R) 
= 2 \sum_{i}\,\lambda_i\,e_i(f) e_i (R)\] and $e_i (R) =  < \nabla R, e_i> = 2 \ric (\nabla f, e_i) = 2 \lambda_i e_i (f)$. So for each $i$, $\,e_i(f) e_i (R) \geq 0$.  Hence $|\nabla R|^2 \leq 2 \Lambda <\nabla f, \nabla R >$. Finally, we recall that $< Z, \nabla f > \geq 0$ (Remark \ref{z-f}). We thus conclude,  from (\ref{eq:f-nabla-f}), that 
 \begin{eqnarray}
F \, \leq\,  \frac{1}{2} <\nabla f, \nabla R>  F_1,  \label{eq:bu2}
\end{eqnarray}
where 
\begin{eqnarray*}
F_1 = (\phi - \phi^{\prime}) (f) + \Lambda  ( \phi + \phi^{\prime})(f) + 4 \Lambda ( \phi + \phi^ {\prime}) - \frac{1}{2} f  (\phi^{\prime\prime} + \phi^{\prime})(f).
\end{eqnarray*}
It follows from (\ref{eq:bu1}) and (\ref{eq:bu2}) that 
\begin{eqnarray}
\int_{M}\, \phi (f) Q  \, \leq \, \frac{1}{2} \int_{M} \, <\nabla f, \nabla R>\,F_1.   \label{eq:bu3}
\end{eqnarray}
A direct computation leads to
\begin{eqnarray*}
&&F_1 = (\phi - \phi^{\prime}) (t) + \Lambda  ( \phi + \phi^{\prime})(t) + 4 \Lambda ( \phi + \phi^ {\prime})(t) - \frac{1}{2} t (\phi^{\prime\prime} + \phi^{\prime})(t)\\
&&=  -\frac{1}{2} \delta (1+\delta)(t-T)^{k+1} e^{\delta t}  - [ \frac{1}{2} (1+2\delta) k - 5(1+\delta) \Lambda -1  + \frac{T-2}{2}\delta ](t-T)^{k} e^{\delta t}\\
&& \hspace{1cm} - k [ \frac{1}{2} (k-1) - 5 \Lambda + \frac{1}{2} T +1] (t-T)^{k-1}e^{\delta t}- \frac{1}{2} T \phi^{\prime\prime}.
\end{eqnarray*}
If we choose $k > 10 \Lambda +2$, the above expression will clearly be negative for $t >T$.  
We have therefore shown  that  $F_1 \leq 0 $ everywhere and $F_1 < 0$ where $f > T$.  Since $ Q \geq 0 $ (Remark \ref{Q+}) and $ < \nabla f, \nabla R> = 2 \ric(\nabla f, \nabla f) \geq 0$ (Lemma \ref{basic}), we conclude from (\ref{eq:bu3}) that  $ < \nabla f, \nabla R> =0$ in the region $ \{ f  > T\}$. But as we have noted earlier in the proof, $|\nabla R|^2 \leq 2 \Lambda <\nabla f, \nabla R >$. Hence $\nabla R =0$ in the region $ \{ f  > T\}$. The analyticity of metric  (\cite{B}, \cite{K}) then implies that $R$ is a constant function. Theorem \ref{main0.5} then follows from \cite{PW1}. 

 \hfill{} \qedbox

\section{Proof of Theorem \ref{main1}}\[\]  \vspace{-1.2cm}

We first show that the Ricci tensor has a zero eigenvalue at any point $p$ in $C$,  then show that the soliton splits in a neighborhood of $p$, which, in turn, implies that the scalar curvature is a constant. 
 
Let $C$ be the critical manifold of minima of $f$. Since $C$ is assumed to be non-degenerate, the Bott-Morse Lemma implies that for any point $p \in C$, there exists an open neighborhood $U$ of $p$ and a diffeomorphism $\phi : U \rightarrow \mathbb{R}^n$ such that $ \phi (U \cap C) = \{ (0, ..., 0,  x_{m+1}, ..., x_n)\}$ , $\phi(p) = 0$ and  $f\circ \phi^{-1} (x_1, ..., x_n) = c+ \frac{1}{4} (x_{1}^2 + ... + x_m^2)$.  

In what follows in this section, unless specified otherwise, the range for the greek leters $\alpha, \beta, ... $  is $1$ to $m$ while that for the latin letters $i, j, ...$ is $m+1$ to $n$. 

We observe that we may assume that for all $\alpha$ and $i$, $g^{\,\alpha i } (p)=0$ . In fact, by making a change of variables, $y_{\alpha} = x_{\alpha}$ and $y_{i} = x_{i}- \sum_{\beta=1}^{m}\, g^{i\, \beta}(p) x_{\beta}$, we see that
in the new coordinates, at $p$, 
$g^{ \alpha \, i } = <\nabla y_{\alpha}, \nabla y_{i}> = 0$ for $\alpha$ and $  i $. Moreover, $ f(y_1, ..., y_m, y_{m+1}, ..., y_n) = c + \frac{1}{4} (y_1^2 +... + y_m^2)$. From now on, we assume in the original coordinates $(x_1, ..., x_n)$, $g^{\alpha i }(p) =0 $ for all $ \alpha $ and $i$. As a consequence, we also have $g_{\alpha i}(p) =0$.  

Next lemma computes the Ricci tensor at $p$.

\begin{lem} \label{ricci-tensor} At $p$, we have 
 $ \ric(p) (\frac{\partial}{\partial x_{\alpha}},  \frac{\partial}{\partial x_{\beta}} )= \frac{1}{2} ( g_{\alpha \beta}(p) - \delta_{\alpha\beta}) $; $ \ric(p) (\frac{\partial}{\partial x_{i}},  \frac{\partial}{\partial x_{j}} )= \frac{1}{2} g_{ij}$; and
$ \ric (p) (\frac{\partial}{\partial x_{\alpha}},  \frac{\partial}{\partial x_{i}} )= 0.$
  \end{lem}
 
\proof Since $ \nabla f =   \frac{1}{2} g^{\alpha \beta } x_\alpha \frac{\partial}{\partial x_{\beta}} + \frac{1}{2} g^{\alpha i } x_{\alpha} \frac{\partial}{\partial x_i}$,  we have at $p$, $ \hess (f) (p) (\frac{\partial}{\partial x_{\alpha}},  \frac{\partial}{\partial x_{\beta}} )= \frac{1}{2} \delta_{\alpha\beta}$, and $ \hess (f) (p) (\frac{\partial}{\partial x_{\alpha}},  \frac{\partial}{\partial x_{i}} )= \hess (f) (p) (\frac{\partial}{\partial x_{i}},  \frac{\partial}{\partial x_{j}} )= 0$. The lemma follows from the soliton equation. 

\hfill{} \qedbox

Let $\mu_{\gamma}^{-1} $ ( $\gamma =1, ..., m$) denote the eigenvalues of the positive definite symmetric matrix $g_{\alpha \beta}(p)$. Then there exists $ (v_{1 \gamma}, \ldots, v_{m \gamma }) \neq 0$ such that  $\sum_{\beta} g_{\alpha \beta} (p) v_{\beta\gamma} = \mu_{\gamma}^{-1} v_{\alpha\gamma}$. Let $v_{\gamma} = \sum_{\alpha} v_{\alpha \gamma} \frac{\partial}{\partial x_{\alpha}}$. The first part of Lemma \ref{ricci-tensor} implies that 
 \begin{eqnarray*}
\ric(p) (v_{\gamma}, v_{\gamma}) &=& \sum_{\alpha, \beta} v_{\alpha \gamma} v_{\beta \gamma} \ric(p) (\frac{\partial}{\partial x_{\alpha}},  \frac{\partial}{\partial x_{\beta}} ) \\
&=& \frac{1}{2} (\mu_{\gamma}^{-1}-1) \sum_{\alpha} (v_{\alpha  \gamma})^2 \\
&=& \frac{1}{2} (\mu_{\gamma}^{-1}-1) \mu_{\gamma} g(p)(v_{\gamma}, v_{\gamma})\\
&=& \frac{1}{2} (1-\mu_{\gamma})  g(p)(v_{\gamma}, v_{\gamma}). 
 \end{eqnarray*}
We conclude from this and the rest of Lemma \ref{ricci-tensor} that the eigenvalues of the Ricci tensor  at $p$ are $\frac{1-\mu_{\alpha}}{2}$, $\alpha=1, ..., m$ and $\frac{1}{2}$ with multiplicity $n-m$.  Since the Ricci tensor is assumed to be semi-positive definite, $\mu_{\alpha} \leq 1$ for each $\alpha$. Of course $\mu_{\alpha} > 0$.  Our goal is to show that $\mu_{\alpha}=1$. 

Now assume $\{e_1, ..., e_{n}\}$ is an orthonormal basis  in a neighborhood of a fixed point $p \in C$ with $D_{e_i}e_{j}\,(p)=0$ for $ 1 \leq i, j, \leq n$. We may assume that each $e_\alpha$ is an eigenvector 
of $\ric$ at $p$ corresponding to the eigenvalue $\frac{1-\mu_{\alpha}}{2}$ for $1\leq \alpha \leq m$ and $e_i$ an eigenvector 
corresponding to $\frac{1}{2}$ for $m+1\leq i \leq n$. 

By our assumption, $D\ric =D^2\ric =0$ at $p$. Hence, for each $ 1 \leq s \leq n$, in the neighborhood of $p$, 
\[ \ric(e_s, e_s) = r_s + \sum_{i,j,k=1}^n  r_{sijk}x_ix_j x_k + h.o.\]
where $r_s$  and $r_{sijk}$ are constants. We make the following observation.

\begin{lem} \label{ricci-tensor2} we have 
 \begin{eqnarray*}
&& r_{\alpha} = \frac{1-\mu_{\alpha}}{2}, \,\,\, \alpha=1, ..., m;\,\,\,\, r_{i} = \frac{1}{2}, \,\,\, i=m+1, ..., n \\
&& \hspace{2cm}   \sum_{\alpha=1}^{m}K_{s \alpha } \mu_{\alpha} = 0,
 \end{eqnarray*}
 where $K_{s\alpha}$ is the sectional curvature of the section spanned by $e_s$ and $e_{\alpha}$.
  \end{lem}

\proof We only need to prove the second line. At $p$, 
\[
(\Delta \ric) (e_s, e_s) = \Delta [\ric(e_s, e_s)] =  0. \]
On the other hand, we have $\Delta \ric = D_{\nabla f} \ric + \ric -2 \sum_{l=1}^n\, Rm(\cdot, e_l, \ric(e_l), \cdot) $ (Lemma 2.1 in \cite{PW2}, see also the proof of Lemma \ref{delta-ricci-square}). Hence,
\begin{eqnarray*}
0 &=& \ric(e_s, e_s)-2 \sum_{l=1}^n\, Rm(e_s, e_l, \ric(e_l), e_s)\\
&=& r_s - 2 \sum_{\alpha=1}^m\, Rm(e_s, e_{\alpha}, \ric(e_{\alpha}), e_s) - 2 \sum_{i=m+1}^n\, Rm(e_s, e_{i}, \ric(e_{i}), e_s) \\
&=& r_s - \sum_{\alpha=1}^m\, (1- \mu_{\alpha}) Rm(e_s, e_{\alpha}, e_{\alpha}, e_s) - \sum_{i=m+1}^n\, Rm(e_s, e_{i}, e_{i}, e_s) \\
&=& \sum_{\alpha=1}^{m}K_{s \alpha } \mu_{\alpha}.
\end{eqnarray*}
\hfill{} \qedbox

\noindent We are now in position to prove Theorem \ref{main1}. 
\vspace{.2in}

{\em Proof of Theorem \ref{main1}} \,\,\,\, It follows from Lemma \ref{ricci-tensor2} and the assumption of nonnegative sectional curvature that $K_{s\alpha} (p) =0$ for all $ 1\leq s \leq n$. So,
 $\ric$(p) vanishes on the subspace spanned by $\{\frac{\partial}{\partial x_{\alpha}} \, | \, \alpha =1, ..., m\}$.  

We first prove that a neighborhood of $p$ splits isometrically as $U \times V$, where $U$ is of  at least $m$ dimensional and $\ric \equiv 0$ on $U$. We have shown that  that $\ric_{\alpha \beta} (p) =0$.
The rest of the argument are along the lines of the proof of Lemma 8.2 in ~\cite{H3} and that of Corollary 2.1 in ~\cite{NT}. Denote by $K(x, t)$ the null space of $\ric(x,t)$, i.e.
\[ K(x,t) = \{w \in T_x\, M \, | \, \ric(x,t)(w) = 0 \} \]
Let $w_0 \in K(p, -1)$ and $\gamma(s)$ a smooth curve starting from $p$. Parallel translating $w_0$ along $\gamma$ gives a vector field $w$ along $\gamma$. Denote  the extension of $w$ to a neighborhood of $\gamma$ still by $w$. Now we project $w$ onto $K(x,t)$ to get a vector field $ v(x,t)$. Then $v(\gamma(s), t) \in K(\gamma(s), t)$. We first show that $D_{\gamma^{\prime}}\, v$ is also in $K(\gamma(s),t)$. We fix 
 an orthonormal 
 basis in $g(t)$, $\{e_1, ..., e_{n}\}$, in a neighborhood of a fixed point $\gamma(s)$ and assume that $e_{i}(\gamma(s))$ are the eigenvectors of $\ric$. For simplicity of notations, we denote $e_i(\gamma(s))$ by $e_i (s)$. 
Since $\ric(v)=0$, $[\frac{\partial}{\partial t}\,\ric]\,(v,v)\,=0$. The evolution equation for Ricci tensor then implies that at $\gamma(s)$, 
\[ ( \Delta \ric ) (v, v) - 2 <\ric (v), \ric(v) >  + 2 \ric (e_i, e_i)\,K(e_i, v) =0 , \]
where the repeated indices are being summed over. Since the sectional curvature $K(e_i, v) \geq 0$ and since $\ric(v)=0$, we deduce that $ (\Delta \ric ) (v, v)  \leq 0 $. Direct computations give 
\begin{eqnarray*}
&&(\Delta \ric ) (v,v) = \Delta [\ric(v,v)] - 4 e_i [ \ric (v, D_{e_i}\,v)] + 2 \ric (v, D_{e_i}D_{e_i}\,v) \\
&& \hspace{5cm} + 2 \ric(v, D_{D_{e_i}\,e_i}\,v) + 2\ric(D_{e_i}\, v, D_{e_i}\,v).
\end{eqnarray*}
Using $ (\Delta \ric ) (v, v) \leq 0 $ and $ \ric (v) =0$, we obtain $\ric(D_{e_i}\, v, D_{e_i}\,v) \leq 0$. Since $\ric$ is positive semi-definite, we conclude that $\ric(D_{e_i}\,v) = 0$ for each $i$, and hence $D_{\gamma^{\prime}}\,v\, \in K(\gamma(s),t)$. As in the proof of Corollary 2.1 in ~\cite{NT}, we conclude that $w \in K(x,t)$.  Since parallel translation preserves inner product, for each fixed $t$, the dimension of $K(x,t)$ is independent of $x$. We then use De Rham's decomposition theorem to conclude that a neighborhood of $p$ splits. 

Note  that $|\nabla f|^2 \geq  f$ on $U \times V$.  In fact, for any $q\in V$, the restriction of $g$ and $f$ on $U \times \{q\}$ gives a soliton on $U \times \{q\}$ with zero $\ric$ tensor. Lemma \ref{basic}(2)  implies that  $ |\nabla_{U\times \{q\}}\,f|^2 = f |_{U\times \{q\}} $, where $\nabla_{U \times \{q\}}\,f$ is the gradient of $f |_{U\times \{q\}}$ with respect to the metric $g |_{U \times \{q\}}$. Since $|\nabla f|^2 \geq |\nabla_{U \times \{q\} }\,f|^2 $, we infer that $|\nabla f|^2 (x, q) \geq  f (x, q)$ for all $ x\in U$. Since $q$ is an arbitrary point in $V$, it follows that 
$|\nabla f|^2 \geq  f$ on $U \times V$.

We now prove that $|\nabla f|^2 \leq  f$  on $U \times V$. Given any point $y \in U\times V$, denote by $\gamma(s)$ the integral curve of $\frac{\nabla f}{|\nabla f|^2}$ such that $\gamma(0) =  y$. Then $f(\gamma(s)) = s + f(\gamma(0))$. On the other hand, using Lemma \ref{basic}(1) (2), we have
\begin{eqnarray*}
\frac{d}{ds}\, |\nabla f|^2 (\gamma(s))\,&=&\, \frac{1}{|\nabla f|^2}\, \nabla f ( |\nabla f|^2 ) = \frac{1}{|\nabla f|^2} ( |\nabla f|^2 - <\nabla f, \nabla R>) \\
&= & \frac{1}{|\nabla f|^2} [ |\nabla f|^2 - 2 \ric ( \nabla f, \nabla f)].
\end{eqnarray*}
Since $\ric(\nabla f, \nabla f) \geq 0$, we obtain $\frac{d}{ds}\, |\nabla f|^2 (\gamma(s))\,\leq 1$.  Integrating this inequality  from $-f(\gamma(0))$ to $s$ and noting that $\nabla f (\gamma(s))=0$ at $s=-f(\gamma(0))$ give us the desired inequality $|\nabla f|^2 \leq  f$.

We have thus proved that $|\nabla f|^2 = f$, which, when combined with Lemma \ref{basic}(2), implies  that $R$ is constant in a neighborhood of $p$. Hence $R$ is constant on the entire $M$. The proof of Theorem \ref{main1} is therefore completed. 

 \hfill{} \qedbox

\end{document}